\documentclass[11pt,a4paper]{report}
\usepackage{amsmath,amsfonts,amsthm,amssymb,amscd,hyperref}
\usepackage[french]{babel}
\usepackage[utf8]{inputenc}
\usepackage[T1]{fontenc}

\def\classification#1{\def\@class{#1}}
\classification{\null}

\DeclareFontFamily{OT1}{rsfs}{}
\DeclareFontShape{OT1}{rsfs}{n}{it}{<-> rsfs10}{}
\DeclareMathAlphabet{\mathscr}{OT1}{rsfs}{n}{it}

\DeclareMathOperator{\mo}{\,mod}

\DeclareMathOperator{\GL}{GL}
\DeclareMathOperator{\diam}{diam}

\DeclareMathOperator{\SL}{SL}
\DeclareMathOperator{\SU}{SU}
\DeclareMathOperator{\SO}{SO}

\DeclareMathOperator{\PSL}{PSL}

\DeclareMathOperator{\Sym}{Sym}

\DeclareMathOperator{\Alt}{Alt}
\DeclareMathOperator{\Sp}{Sp}

\DeclareMathOperator{\Gal}{Gal}

\DeclareMathOperator{\Cl}{Cl}

\newtheorem{prop}{Proposition}
\newtheorem{theo}[prop]{Théorème}

\newtheorem{lem}[prop]{Lemma}

\theoremstyle{remark}

\numberwithin{equation}{section}
\title{Groupes, courbes et croissance: synthèse des travaux présentés en vue d'une Habilitation à Diriger des Recherches}
\author{Harald Andrés Helfgott}
\begin{document}
\maketitle
\tableofcontents
\chapter{Introduction}
Le thème central de ma recherche est le dénombrement, la distribution et la croissance
d'objets discrets dans des structures algébriques. Les travaux décrits dans
ce mémoire peuvent être classifiés de la manière suivante:
\begin{enumerate}
\item La croissance dans les groupes: \cite{MR2415382}, \cite{MR2781932}, \cite{GH1}, \cite{GH2},
\cite{HS};
\item Les points entiers et rationnels sur des courbes :\cite{MR2220098}, \cite{MR2508647}, \cite{MR2394896};
\item Combinatoire additive classique: \cite{MR2773330}, \cite{MR2764161};
\item Théorie analytique et probabiliste des nombres: \cite{MR2394896}, \cite{MR2437966};
\item Théorie algorithmique des nombres: \cite{TCH}.
\end{enumerate}
Tous les travaux cités ci-dessus sont postérieurs à la thèse \cite{Hth}. L'article \cite{MR2099831}, issu de cette dernière, constitue un premier pas vers \cite{MR2220098}.

(a) La combinatoire additive, l'expansion dans les groupes arithmétiques et l'étude du diamètre des groupes finis avaient évolué comme trois
domaines disjoints. Leur liaison dans \cite{MR2415382} -- laquelle, en particulier, a établi la conjecture de \textsc{Babai} sur les petits diamètres pour  les groupes $G=\PSL_2(\mathbb{Z}/p\mathbb{Z})$ -- a initié une période d'intense activité \cite{MR2415383}, \cite{BGSU2},  \cite{Din}, \cite{MR2781932}, 
\cite{GH1}, \cite{Var}, \cite{PS}, \cite{BGT}, \cite{GH2}, \cite{SGV}.\footnote{Cette liste ne prétend pas être complète.} Ces travaux, néanmoins, étaient tous restreints au
cas de groupes linéaires de rang borné. Très récemment, \cite{HS} a donné
une borne quasipolynomiale pour le diamètre du groupe symétrique
$\Sym(n)$. Le traitement de ce cas fondamental se base sur un développement d'idées
dans plusieurs domaines, comme la croissance dans les groupes, la structure des groupes de permutation et les marches aléatoires.

(b) Soit $C$ une courbe définie par une équation avec coefficients dans $\mathbb{Q}$. Combien de points entiers et rationnels jusqu'à une certaine
hauteur peut-il y avoir sur $C$? Il y a des bornes \cite{MR1016893} prouvées à partir de la géométrie de la courbe comme sous-ensemble de $\mathbb{R}^2$, 
et des bornes (\cite{MR0186624}, \cite{MR895285}, \cite{MR1328329}, \cite{MR2099831}) 
qui exploitent la géométrie, voire la norme naturelle, du réseau de Mordell-Weil de $C$ (en genre positif). La méthode suivie dans \cite{MR2220098} (avec \textsc{A. Venkatesh}) réussit
à combiner des idées des deux types. Parmi les applications se trouvent des bornes non-triviales pour le nombre de courbes elliptiques
de conducteur donné et pour la 3-torsion dans des groupes de classe \cite{MR2220098}. Une approche similaire est un des éléments qui conduisent
à des résultats sur les valeurs $f(p)$ ($f$ cubique) sans facteurs
carrés (\cite{MR2394896}, \cite{MR2437966}; voir ci-dessous). Par ailleurs,
le transfert d'idées sur la géometrie des courbes dans $\mathbb{R}^2$ dans le contexte du {\em larger sieve} de \textsc{Gallagher} \cite{MR0291120}
conduit à une forte dichotomie entre cardinalité et structure algébrique en $\mathbb{Z}^2$ \cite{MR2508647}.

(c) Dans le cas classique, c'est-à-dire les groupes abéliens, la combinatoire additive peut arriver à combiner des idées de l'analyse de Fourier avec la
théorie analytique des nombres. Une telle combinaison a conduit à une amélioration sensible \cite{MR2773330} du résultat de \textsc{Green} \cite{MR2180408} sur
les progressions arithmétiques de longueur $3$ dans tout sous-ensemble dense des  nombres premiers. Une approche arithmétique-combinatoire s'est montrée suffisante pour donner des bornes explicites pour un théorème d'incidence sur $\mathbb{F}_p$ \cite{MR2764161} . 

(d) Les articles \cite{MR2394896}, \cite{MR2437966} prouvent une conjecture  
d'Erd\H{o}s \cite{MR0056635} selon laquelle, pour tout polynôme cubique $f\in \mathbb{Z}\lbrack x\rbrack$ satisfaisant des conditions locales nécessaires, il y a un nombre infini de nombres premiers $p$ tels que $f(p)$ soit sans facteurs carrés.  La preuve repose en partie sur
des méthodes de comptage de points sur des courbes (similaires à celles dans \cite{MR2220098}) et en partie
sur des méthodes probabilistes d'un type peut-être nouveau en théorie des nombres.
Par exemple, \cite{MR2394896} utilise le fait qu'il y a deux distributions incompatibles pour deux
valeurs reliées -- l'argument et la valeur de $f(x)$. Dans un tel cadre, la preuve dans \cite{MR2437966} utilise une conséquence de la modularité de $y^2 = f(x)$: le terme d'erreur dans 
l'estimation de \textsc{Hasse}-\textsc{Weil} existe réellement
en moyenne, et donc peut être exploité.

(e) Il est possible de vérifier rapidement si un nombre $n$ est premier -- 
ceci est classique si on permet une méthode probabiliste, et connu depuis 
\cite{MR2123939} de manière déterministe. Par contre, s'il s'agit de 
{\em trouver} un nombre premier entre $N$ et $2N$, il reste une grande 
différence entre les deux régimes: si on permet un algorithme probabiliste,
le problème se réduit a celui de la 
vérification en temps $O(\log n)$ de manière triviale, mais, même sous l'hypothèse de Riemann, il n'est pas évident
de faire mieux que $O\left(n^{1/2 + \epsilon}\right)$ de manière
déterministe. Dans \cite{TCH} -- un article issu du projet {\em Polymath} de collaboration massive
-- une méthode déterministe basée sur une 
analogie avec des approches combinatoires au problème des diviseurs de
\textsc{Dirichlet} (voire le problème du cercle de \textsc{Gauss})
accomplit la tâche en temps $O\left(n^{1/3+\epsilon}\right)$ sans utiliser
l'hypothèse de Riemann, en supposant que le nombre de nombres premiers
entre $N$ et $2N$ soit impair.

{\bf Note (12 déc 2011).} Je remercie les membres du jury -- N. Anantharaman, E. Fouvry,
M. Hindry, E. Kowalski et P. Pansu -- pour ses questions et ses commentaires.
\chapter{Croissance et diamètre}

\section{Résultats principaux}

Soit $G$ un groupe, $A\subset G$ un sous-ensemble fini. On écrit
\[\begin{aligned}
A^{-1} &= \{g^{-1} : g\in A\},\;\;\;\;\;\;\;\;\;\;\;\;\;\;\;\;\; \;\;\;\;\;
A\cdot A = \{g_1\cdot g_2: g_1,g_2 \in A\},\\ 
A\cdot A\cdot A &= \{g_1\cdot g_2\cdot g_3: 
g_1, g_2, g_3 \in A\},\;\;\;\;
A^k = \{g_1 g_2 \dotsb g_k : g_1, g_2,\dotsc,g_k \in A\}.\end{aligned}\]
Rappelons que $A$ est un {\em ensemble de générateurs}
de $G$ (ou: $\langle A\rangle = G$) 
si tout élément de $G$ peut être écrit comme un produit d'éléments de $A$ et
ses inverses, i.e., si  $G = \bigcup_{k\geq 0} (A \cup A^{-1})^k$. A partir de maintenant nous
supposons, pour simplifier, que $A = A^{-1}$ et $e\in A$, où $e$ est l'identité de $G$.

Si $\langle A\rangle = G$, on définit le {\em graphe de Cayley} (non-orienté) $\Gamma(G,A)$,
dont l'ensemble de sommets est $G$ et l'ensemble d'arêtes est $\{\{g, a g\}: g\in G, a\in A\}$.
Cette définition permet de formuler certains concepts de manière géométrique.
La longueur d'un chemin dans un graphe est tout simplement le nombre d'arêtes dans le chemin;
la {\em distance} $d(x,y)$ entre deux points est le longueur du chemin le plus court entre eux. 
Le {\em diamètre} $\diam(\Gamma)$ d'un graphe $\Gamma$ est le maximum de $d(x,y)$, où $x$, $y$ parcourent tous les sommets du graphe. Il est facile
de voir que la distance entre $g\in G$ et $h\in G$ dans un graphe de Cayley est l'entier
non-négatif $k$ 
minimal tel que $g h^{-1}\in A^k$. Donc, le diamètre de $\Gamma(G,A)$ est l'entier non-négatif
$k$ minimal tel que $G = A^k$.

Le diamètre $\diam(G)$ d'un groupe fini $G$ est le maximum de $\diam(G,A)$ quand $A$ parcourt tous les ensembles de générateurs $A$ de $G$.

\begin{theo}[\cite{MR2415382}]\label{theo:carr}
Soit $G = \SL_2(\mathbb{Z}/p\mathbb{Z})$. Alors
\[\diam(G) \leq (\log |G|)^{O(1)},\]
où la constante implicite est absolue.
\end{theo}
(Ici et dans le reste, $|S|$ dénote le nombre d'éléments d'un ensemble $S$.) La preuve de ce théorème repose sur la proposition suivante:
\begin{prop}$($\cite["Key Proposition"]{MR2415382}$)$\label{prop:key}
Soit $G = \SL_2(\mathbb{Z}/p\mathbb{Z})$. Soit $A\subset G$,
$\langle A\rangle = G$.
Si $|A|<|G|^{1-\epsilon}$ pour $\epsilon>0$, alors \[|A\cdot A\cdot A| \geq |A|^{1+\delta},\]
où $\delta>0$ dépend seulement de $\epsilon$.
\end{prop}

L'itération de cette proposition conduit directement au Théorème \ref{theo:carr}.
La Proposition \ref{prop:key} a été vite utilisée pour construire des graphes expanseurs 
(\cite{MR2415383}; voir aussi les applications dans \cite{MR2587341}) et généralisée à des
autres groupes du même type de Lie ($\SL_2(\mathbb{F}_q)$ \cite{Din} et $\SU_2$ \cite{BGSU2}). Dans \cite{MR2781932}, j'ai généralisé Prop.\ \ref{prop:key} et 
Théo.\ \ref{theo:carr} à $\SL_3(\mathbb{Z}/p\mathbb{Z})$, tout en traçant les lignes à suivre pour
prouver un énoncé plus général. Il y a eu une série de généralisations ultérieures - \cite{GH1} et en suite  \cite{BGT} et \cite{PS} - de telle manière que maintenant on connait Prop.\ 
\ref{prop:key} et  Théo.\ \ref{theo:carr} pour tous les groupes simples de
type de Lie et rang borné.
(Autrement dit - le rang peut être arbitrairement grand, mais \cite{BGT} et \cite{PS}, même
comme \cite{GH1}, gardent une forte dépendance du rang dans ses bornes.)
Des résultats similaires ont été prouvés pour les groupes résolubles \cite{GH2}. 

En général, un énoncé comme le Théorème \ref{theo:carr} doit être vrai pour tous les groupes 
simples finis et non-abéliens (conjecture de \textsc{Babai}). Il restait donc à prouver deux
cas:  les groupes simples de type de Lie et rang non-borné, et les groupes alternés $\Alt(n)$.
Les deux problèmes sont reliés: $\Alt(n)$ apparaît dans les groupes de Weil des groupes
de type de Lie, et la généralisation la plus naïve de la proposition \ref{prop:key} est fausse
dans les deux cas; il y a même des contre-exemples similaires dans un cas et l'autre 
\cite{PS}, \cite{PPSS}. 

Le résultat suivant est très récent. Sa preuve présente une série de
nouveaux arguments utilisant à la fois des idées inspirées de \cite{MR2781932} comme des résultats antérieurs dans le
domaine des groupes de permutation.
\begin{theo}[\cite{HS}]\label{theo:curr}
Soit $G = \Alt(n)$ or $G=\Sym(n)$. Alors
\[\diam(G) \leq e^{(\log n)^{O(1)}} = e^{(\log \log |G|)^{O(1)}}.\]
 Plus précisément, 
$\diam(G)\leq \exp(O((\log n)^4 (\log \log n)))$,
où la constante implicite est absolue.
\end{theo}
 Par \cite{BS92}, la borne
\[\diam(G)\leq \exp(O((\log n)^4 (\log \log n)))\] suit immédiatement aussi pour tout groupe
transitif de permutation sur $n$ éléments.

\section{Domaines et Perspectives}
Jusqu'à \cite{MR2415382}, chacun des trois domaines suivants paraissait disjoint des autres. Jusqu'à \cite{HS}, le dernier l'était encore.

\begin{enumerate}
\item\label{it:diasp}
 {\em Diamètre et spectre dans les groupes arithmétiques.} Il y avait quelques résultats
classiques pour des générateurs spécifiques.  Par exemple, soit
\begin{equation}\label{eq:knop}
A = \left\{ \left(\begin{array}{cc} 1 &1\\0 &1\end{array}\right),
\left(\begin{array}{cc} 1 &0\\1 &1\end{array}\right)
\right\} .\end{equation}
Le théorème de \textsc{Selberg} sur le trou spectral en $\SL_2(\mathbb{Z})\backslash \mathbb{H}$ (\cite{MR0182610}) implique que $\{\Gamma(\SL_2(\mathbb{Z}/p \mathbb{Z}), A)\}_{p\geq 5}$ 
est une famille de graphes expanseurs (voir, e.g., \cite{MR1308046}, Thm.\ 4.4.2, (i)). En particulier,
$\diam(\Gamma(\SL_2(\mathbb{Z}/p\mathbb{Z},A)))\ll \log p$. Malheureusement, cet argument
ne marche pas pour $A$ plus général. Par exemple, il n'y avait pas des bonnes bornes pour
$\diam(\Gamma(\SL_2(\mathbb{Z}/p\mathbb{Z},A)))$ avec, par exemple,
\begin{equation}
A = \left\{ \left(\begin{array}{cc} 1 &3\\0 &1\end{array}\right),
\left(\begin{array}{cc} 1 &0\\3 &1\end{array}\right)
\right\} .\end{equation}
(un exemple favori de \textsc{Lubotzky}).

\item\label{it:combadd} {\em Combinatoire additive.} Soit $A$ un sous-ensemble des entiers (ou d'un groupe
abélien). Que pouvons-nous dire sur la taille de $A+A=\{x+y: x,y\in A\}$? Si $A+A$ n'est
pas beaucoup plus grand que $A$, que pouvons-nous dire sur la structure de $A$?
(\textsc{Freiman}, \textsc{Ruzsa},\dots)

\item\label{it:gracay} {\em Graphes de Cayley des groupes de permutation.} Soit $G<\Sym(n)$ un groupe de permutation,
$A$ un ensemble de générateurs. Le diamètre $\diam(\Gamma(G,A))$ est-il polynomial
($= n^{O(1)}$)? Quasipolynomial ($= \exp\left((\log n)^{O(1)}\right)$) (\cite[Conjecture 1.6]{BS92})? Il avait des bornes
sous certaines conditions \cite{BBS04}, des bornes pour des générateurs aléatoires
\cite{MR1208801}, \cite{MR2298365}, et une réduction au cas $G=\Sym(n)$ 
\cite{BS92}, mais les bornes pour
des générateurs arbitraires restaient très loin des conjectures (\cite{BS88}, borne exponentielle en 
$\sqrt{n}$).
\end{enumerate}

L'article \cite{MR2415382} a résolu un problème de type (\ref{it:diasp}) 
pour $A$ complètement arbitraire en utilisant des méthodes de type
(\ref{it:combadd}). Ceci a mené à des développements sur le trou spectral (la plus petit valeur propre non-nulle du Laplacien du graphe de Cayley): \cite{MR2415383}, \cite{MR2587341}, etc. 
Récemment, \textsc{Bourgain}, \textsc{Gamburd} et \textsc{Sarnak} \cite{BGSup}
ont réussi à inverser la procédure utilisée jadis pour (\ref{it:diasp}): ils déduisent
des trous spectraux sur les surfaces $\Gamma\setminus \mathbb{H}$ à partir
du 
résultat combinatoire dans \cite{MR2415382} (via \cite{MR2415383}).

La preuve  dans \cite{HS} lie les méthodes en (\ref{it:gracay}) avec les arguments développés à
partir de \cite{MR2415382} et \cite{MR2781932}. Une idée motrice est
l'examen de la croissance comme un concept relatif, c'est-à-dire, en termes
d'une action $G\to X$, pas seulement comme un phénomène dans un groupe $G$. Les arguments en
 \cite{MR2415382}, \cite{MR2781932}, \cite{Din}, \cite{GH1},
 \cite{BGT}, \cite{PS} peuvent être encadrés (à partir de \cite[\S 4]{MR2415382})
 en termes de l'action de $G$ sur lui-même par conjugaison ou (dans \cite[\S 4.1 et \S 7.1]{MR2781932}) simplement par multiplication; en rétrospective,
\cite{BBS04} et autres travaux dans le domaine (\ref{it:gracay}) peuvent
être mis dans un cadre similaire en relation à l'action naturelle de $G=\Sym(n)$ sur $X=\{1,2,\dotsc,n\}$. Les deux types d'action
se mêlent dans \cite{HS}.

\section{Outils et idées}
\subsection{Analogies. Adaptations.} Deux analogies ont joué un rôle important dans le
développement  
de ce qu'on peut nommer la combinatoire arithmétique non-commutative -- spécialité qui a ses origines dans \cite{MR2415382} et \cite{MR2501249}.
\begin{enumerate}
\item {\em Analogues non-commutatifs d'énoncés de combinatoire additive.}
Pour prouver la Proposition \ref{prop:key}, on prouve $|A^k|\geq |A|^{1+\delta'}$, $k$ une
constante, et on en déduit que $|A \cdot A\cdot A|\geq |A|^{1+\delta}$; ce dernier pas est
possible grâce à l'inégalité du triangle de \textsc{Ruzsa}, généralisée
aux groupes non commutatifs dans \cite[\S 2.3]{MR2415382} et \cite[\S 3]{MR2501249}. L'article \cite{MR2501249}
contient aussi une généralisation aux groupes non-commutatifs du théorème de 
\textsc{Balog}-\textsc{Szemerédi}-\textsc{Gowers} (voir aussi \cite{MR2155059}); 
cette généralisation a été utilisée ensuite 
dans le passage d'ensembles aux mesures
dans \cite[\S 3]{MR2415383}.

\item\label{it:anadeux} {\em Analogues pour des ensembles d'énoncés sur des sous-groupes.}
Un bon nombre de résultats sur des sous-groupes d'un groupe $G$ peuvent être adaptés à des sous-ensembles $A$ de $G$. Un exemple simple mais central est donné par le théorème
d'orbite-stabilisateur.
\begin{lem}[Théorème orbite-stabilisateur pour ensembles]\label{lem:kot}
Soit $G$ un groupe qui agit sur un ensemble $X$. Soit $x\in X$, et soit $A\subset G$ non-vide.
Alors
\[|A A^{-1} \cap G_x| \geq \frac{|A|}{|x^A|}.\]
De plus, pour tout $B\subset G$,
\[|A B| \geq |A\cap G_x| |x^B|.\]
\end{lem}
Ici $G_x$ denote le stabilisateur de $x$ dans G, et $x^B$ denote l'orbite de $x$ sous $B\subset G$.
Le théorème d'orbite-stabilisateur usuel est le cas spécial $A=B=H$, $H$ un sous-groupe 
de $G$. La preuve, à la fois ici comme dans le cas traditionnel, est un simple exercice.

Ce principe est déjà plus ou moins implicite dans \cite[Prop.\ 4.1]{MR2415382}, dans le sens que la tension entre le centralisateur $C(g)$ (le stabilisateur d'un élément $g\in G$ sous l'action de
$G$ sur soi-même par conjugaison) et le nombre de classes de conjugaison (le nombre
d'orbites dans la même action) est décrite et exploitée.

L'analogie (\ref{it:anadeux}) a été appliquée aussi dans \cite{BG}, 
\cite{Hrushovski} (théorie des modèles) et \cite{BGT}, entre autres. Deux
régimes sont utilisés couramment: celui des ensembles $A$ avec $|A A A|\leq |A|^{1+\delta}$
(dans \cite{MR2415382}, \cite{MR2781932}, \cite{GH1}, etc.) et celui des
ensembles appelés {\em approximate groups} (\cite{MR2501249}, \cite{BG}, \cite{BGT}, etc.). Les deux régimes sont à peu près
équivalents (\cite[Thm. 3.8]{MR2501249}).
\end{enumerate}

La preuve de \cite{HS} confirme l'importance de l'analogie (\ref{it:anadeux}) et suggère qu'elle dépasse les deux cadres utilisés jusqu'à maintenant. Considérons, en particulier, la généralisation
dans \cite[\S 5]{HS} du {\em splitting lemma} de \textsc{Babai} \cite[\S 3]{Bab82}.
\begin{itemize}
\item Il était déjà clair que des assertions sur des sous-groupes basées sur des arguments
quantitatifs robustes pouvaient être adaptées a des sous-ensembles (voir Lemma \ref{lem:kot}).
Ce qui se révèle robuste dans \cite[\S 5]{HS}, par contre, n'est pas un argument
quantitatif, mais un argument probabiliste, grâce au fait qu'une distribution uniforme puisse être
simulée par un processus de complexité bornée, à savoir, une marche aléatoire. 
\item Les énoncés dans \cite[\S 5]{HS} ne paraissent pas entrer de façon naturelle ni dans le formalisme $|A A A|\leq |A|^{1+\delta}$ ni dans celui des ``approximate groups''.
Il s'agit plutôt (comme dans \cite[\S 4]{MR2415382}, \cite[\S 5]{MR2781932}
et ailleurs) d'énoncés qui remplacent parfois $A$ par $A^k$, là 
 où un énoncé sur des sous-groupes $H$ utiliserait toujours $H$
(car $H^k = H$). Ce cadre paraît être le plus simple et flexible.
\end{itemize}

\subsection{Intersections avec des sous-varietés}
Pour $G$ un groupe algébrique linéaire défini sur un corps $K$, un {\em tore} est un groupe
algébrique isomorphe au produit d'un nombre fini de copies du groupe
multiplicatif $\mathbb{G}_m$. Un élément $g\in G(\overline{K})$ est
régulier 
semi-simple si et seulement si son centralisateur est un tore maximal.
(Par exemple, si $G = \SL_n$, un tore est un sous-groupe 
isomorphe (par conjugaison dans
$G(\overline{K})$) au groupe des matrices diagonalisables, et
un élément est régulier semi-simple si et seulement si il a $n$ valeurs
propres distinctes.)

\begin{prop}\label{prop:zut}$($\cite[Cor. 5.4 et Cor. 5.10]{MR2781932}$)$
Soit $G\subset \GL_n$ un groupe classique\footnote{Dans ce contexte, $\SL_n$,
$\SO_n$ ou $\Sp_{2n}$. La preuve dans \cite[\S 5]{MR2781932} s'applique aussi
(au moins) à tous les groupes de Chevalley (\textsc{Gill}, non publié).
L'énoncé de (\ref{eq:dusteloj} est prouvé en [\S 5]{MR2781932} seulement
pour $\SL_n$; la même preuve est valide pour tous les groupes classiques.} 
sur un corps $K$. Soit $A\subset G(K)$, $\langle A\rangle = G(K)$, 
$A = A^{-1}$, $e\in A$.

Il y a une constante $k$ qui dépend seulement de $n$ telle que, pour tout 
tore maximal $T/\overline{K}$,
\begin{equation}\label{eq:unustelo}
|A\cap T(K)| \ll_n |A^k|^{\dim(T)/\dim(G)}.
\end{equation}
En plus, il y a un $g\in A^k$ régulier semi-simple tel que
\begin{equation}\label{eq:dusteloj}
|A^k \cap T(K)| \gg_n \frac{|A|}{|A^{k}|} \cdot |A|^{\dim(T)/\dim(G)}
\end{equation}
où $T$ est le centralisateur de $g$.
\end{prop}
Comme il a été remarqué dans \cite[\S 5.3, Remark]{MR2781932}, la preuve en
\cite{MR2781932} donne (\ref{eq:dusteloj}) pas seulement pour un seul $g$,
mais pour beaucoup de $g$ possibles. Dans \cite{MR2781932},
des résultats du type (\ref{eq:unustelo}) ont été aussi prouvé pour
des tores non maximaux et autres sous-groupes. 

\cite{PS} et \cite{BGT} ont amélioré la Proposition \ref{prop:zut}
en prouvant (\ref{eq:dusteloj}) pour {\em tout} élément $g\in A$. Dans le
cas de \cite{BGT}, cette amélioration a été extraite de \cite{MR2813339}.
\footnote{\cite{MR2813339} existait en manuscrit depuis 1998. \textsc{Hrushovski} a joué un rôle
dans sa diffusion.} 
\textsc{Larsen} et \textsc{Pink} avaient prouvé dans \cite{MR2813339}.
des résultats analogues à (\ref{eq:unustelo}) -- pour des sous-groupes $H$
en place de sous-ensembles $A$ arbitraires, mais, par contre, aussi pour des variétés $V$
arbitraires en place de $T$.  En particulier, en prenant
$V=\overline{\Cl(g)}$ (la clôture de Zariski d'une classe de conjugaison)
et utilisant un argument du type orbite-stabilisateur, \textsc{Larsen} et
\textsc{Pink}
 avaient prouvé dans \cite[Thm. 6.2]{MR2813339} une version de
(\ref{eq:dusteloj}) (pour $A$ un sous-groupe) valide pour tout $g\in A$.

Dans \cite{MR2415382}, \cite{MR2781932},
 \cite{GH1}, \cite{BGT}, \cite{PS},
la stratégie générale est la suivante: si $A\subset G$ satisfait $|A\cdot
A\cdot A|\leq |A|^{1+\delta}$, $\delta>0$ petit (c'est-à-dire, si $A$ es
un contre-exemple hypothétique à la Prop.~\ref{prop:key}), alors les bornes
inférieure et supérieure dans la Prop.~\ref{prop:zut} sont très proches
l'une de l'autre. En étudiant des conjugués du tore $T$, on utilise cette
connaissance excessive de la grandeur de $A\cap T(K)$ 
pour arriver à une contradiction.

\subsection{Somme, produit, incidence, pivotement}

Mais comment arriver à une telle contradiction? Dans \cite{MR2415382}, 
une identité pour des traces à été utilisée pour construire un sous-ensemble
de $\mathbb{Z}/p\mathbb{Z}$ qui ne croît ni sous addition ni sous
multiplication par lui-même, ainsi allant à l'encontre du théorème
somme-produit (\cite{MR2053599}, \cite{Koarx}; \S \ref{subs:soprod}).
L'identité est valide seulement pour $\SL_2$ et autres groupes de rang $1$.

Dans \cite[\S 3]{MR2781932}, une idée présente dans des preuves de théorèmes de
type somme-produit (\cite{MR2359478}, \cite{MR2481734}) a été développée et
généralisée à des groupes (abéliens) 
agissant sur des groupes (abéliens ou non-abéliens); \cite{MR2781932}
a nommé cette idée {\em pivoting} (pivotement). Le plan original avait été
de remplacer toute utilisation du théorème somme-produit par
des résultats comme \cite[Prop.\ 3.1]{MR2781932}, qui utilise le pivotement
dans le contexte d'action de groupes. Ce plan n'a pas réussi pour
les groupes simples à
l'époque (en rétrospective, cela est dû au fait que (\ref{eq:dusteloj}) n'était pas
prouvé pour tout $g$), mais il a marché (avec effort) 
pour des groupes résolubles \cite{GH2}. \'A la place, \cite{MR2781932} a
utilisé (et développé) des théorèmes d'incidence (dans le sens
de Szemerédi-Trotter et \cite{MR2053599}). Cette approche s'est avérée 
utile pour $\SL_n$ dans \cite{GH1}; comme il a été remarqué dans \cite{GH1}, elle
est valide aussi pour $\Sp_{2n}$, mais, semble-t-il, pas pour $\SO_n$.

L'argument dans \cite{BGT} est tout à fait basé sur le pivotement. En
\cite{PS}, une méthode essentiellement équivalente a été trouvée par une route
indépendante, basée en partie sur des résultats sur la croissance en 
sous-groupes apparentés à ceux dans \cite[\S 7.1]{MR2781932}. 

Finalement -- un argument crucial dans \cite{HS} (le pas de {\em création},
\cite[Lem. 6.1]{HS}) a une certaine relation avec ce cercle d'idées. Le
contexte relatif est légèrement différent -- cette fois, l'action
est celle d'un normalisateur $N_G(H)$ sur un sous-groupe $H<G$ (par 
conjugaison).

\chapter{Points entiers dans le plan}

\section{Points entiers sur des courbes}\label{sec:courbes}

\subsection{Répulsion dans le réseau de \textsc{Mordell}-\textsc{Weil}}

Soit $C$ une courbe donnée par une équation
\begin{equation}\label{eq:clunt}f(x,y)=0,\end{equation}
où $f\in \mathbb{Q}\lbrack x,y\rbrack$. Le {\em genre} $g$ de $C$ est un
entier $\geq 0$, égal au nombre de trous dans la surface
de Riemann compacte associé.
 Si le genre de $C$ est strictement positif, le nombre de points entiers
sur $C$ (c'est-à-dire, des points $(x,y)\in \mathbb{Z}^2$ satisfaissant (\ref{eq:clunt})) est fini (théorème de \textsc{Siegel}).

Donner une borne supérieure pour le nombre de ces points est une autre
affaire. Les méthodes de \textsc{Siegel} ne sont pas effectives -- elles ne
peuvent donner aucune borne.  Les
travaux de \textsc{Baker} donnent une borne sur la hauteur
(naïve)\footnote{La {\em hauteur naïve} $h_x(P)$ d'un point
entier $P=(x,y)$ sur une courbe dans le plan est tout simplement $\log |x|$.} des points entiers sur $C$,
mais ceci ne donne q'une borne astronomique sur le nombre de points entiers.

L'objet d'intérêt est autant le nombre total de points entiers sur $C$ que
 le nombre de tels
points entiers dans le carré $\lbrack -N,N\rbrack^2$:
\begin{equation}\label{eq:koval}
\{(x,y)\in \lbrack -N,N\rbrack^2 \cap \mathbb{Z}^2 : f(x,y)=0\}.\end{equation}

Une famille de méthodes (\cite{MR1544776}, \cite{MR0337857}, \cite{MR778171}, \cite{MR1016893}) utilise la géométrie de l'ensemble de points $(x,y)\in \mathbb{Q}^2$
vu comme un sous-ensemble de $\mathbb{R}^2$. Les meilleures bornes de ce type sont
celles de \cite{MR1016893} :
\begin{equation}\label{eq:poleda}
|\{(x,y)\in \lbrack -N,N\rbrack^2 \cap \mathbb{Z}^2 : f(x,y)=0\}| \ll_{d,\epsilon} N^{1/d + \epsilon}
\end{equation}
où $d$ est le degré de $C$. (Ici la notation $\ll_{\epsilon,d}$ implique, comme d'habitude, que
la constante implicite dépend seulement de $\epsilon$ et $d$.)

La méthode dans \cite{MR1016893} peut être vue comme un analogue en haut degré du 
principe suivant: si un triangle a des sommets $P$, $Q$, $R$ avec coordonnées entières, alors
son aire appartient à $(1/2) \mathbb{Z}$, et donc elle est au moins $1/2$; par contre, si $P$, $Q$, $R$
sont des points proches les uns aux autres sur une courbe $C$, alors ils sont presque
colinéaires, et donc l'aire du triangle $PQR$ est très petite. Ceci mène a une contradiction si $P$,
$Q$, $R$ sont très proches. Donc il ne peut pas y avoir trop de tels points dans une boîte
$\lbrack -N,N\rbrack^2$.

Une deuxième famille de méthodes (\cite{MR0186624}, \cite{MR895285}, \cite{MR1328329}, \cite{MR2099831}) se base sur la géométrie du réseau de Mordell-Weil - 
c'est-à-dire, le groupe abélien constitué par les points rationnels de la Jacobienne $J_C$ de la courbe elliptique $C$. Dans le cas d'une courbe elliptique, la Jacobienne est isomorphe à $C$ elle-même, et donc on peut travailler directement
avec $C$: il y a une loi d'addition sur l'ensemble
\[C(\mathbb{Q}) = \{(x,y) \in \mathbb{Q}^2 : f(x,y)=0\}\]
et cette loi fait de $C(\mathbb{Q})$ un groupe abélien (\textsc{Poincaré}) avec un nombre fini de générateurs
(\textsc{Mordell}); alors $C(\mathbb{Q}) \sim \mathbb{Z}^r \times T$, où
$T$ est fini. Il y a une {\em hauteur canonique} sur $C(\mathbb{Q})$, c'est-à-dire, une
fonction $\hat{h}:C(\mathbb{Q})\to \mathbb{R}$ telle que
\begin{itemize}
\item la différence entre $\hat{h}(P)$ et la hauteur naïve $h(P) = \log \max(|a|,|b|)$ (pour 
$P = (x,y)$, $x =a/b$, $\gcd(a,b)=1$) est bornée par une constante $O_C(1)$ qui dépend 
seulement de $C$;
\item $\hat{h}$ est une forme quadratique -- en d'autres termes, $\hat{h}(P+Q)+\hat{h}(P-Q)= 2 \hat{h}(P) + 
2 \hat{h}(Q)$;
\item $\hat{h}(P)=0$ si et seulement si $P$ est de torsion, i.e., $P\in T$.
\end{itemize}
Il est donc possible de définir une application $\iota:C(\mathbb{Q})\to \mathbb{R}^r$ avec
$|\iota(P)| = \sqrt{\hat{h}(P)}$. Cette procédure identifie $C(\mathbb{Q})$ modulo $T$ avec un réseau de
$\mathbb{R}^r$.

Les points entiers sur $C$ sont envoyés par $\iota$ sur un sous-ensemble du réseau 
$\iota(C(\mathbb{Q}))$. Quelles sont les propriétés de ce sous-ensemble?

Si $P, Q\in C(\mathbb{Q})$ sont deux points entiers distincts avec $|\iota(P)|\sim (1\pm \epsilon)
|\iota(Q)|$, alors $\iota(P)$ et $\iota(Q)$ doivent se {\em repousser}: l'angle $\measuredangle
\iota(P) O \iota(Q)$ (où $O$ est l'origine) est $\geq 60^\circ -O(\epsilon)$ \cite[Lem.\ 4.17]{MR2099831}. Un phénomène de répulsion a été déjà étudié par Mumford dans le cas de points
rationnels sur des courbes de genre $\geq 2$ \cite{MR0186624}. Le phénomène
observé dans \cite[\S 4]{MR2099831} avait (en essence) été déjà remarqué dans \cite{MR895285}, puis formulé et exploité de manière différente.  Dans \cite{MR2099831} (et \cite{MR2220098}, \cite{MR2394896}), la répulsion est vue de manière géométrique et exploitée par le moyen de bornes sur l'empilement de sphères. (Le nombre de points possibles sur une sphère 
$S_{r-1}\subset \mathbb{R}^r$ séparés par des angles $\geq 60^\circ$ est égal au nombre de
sphères qu'on peut placer autour d'une sphère du même rayon, i.e., le {\em kissing number}
(nombre d'osculation). 
Les meilleures bornes asymptotiques pour le nombre d'osculation, soit pour 
$60^\circ$, soit pour un angle arbitraire, sont celles données par \textsc{Kabatjanski{\u\i}} et
\textsc{Leven{\v{s}}te{\u\i}n}
\cite{MR0514023}. )

Dans \cite{MR2220098}, cette répulsion est dérivée à partir de hauteurs locales $\lambda_v$,
donc de manière plus proche de \cite{MR895285}, \cite{MR1328329}; il s'agit de prouver
des inégalités de type non-archimédien approché pour les hauteurs locales $\lambda_v$
({\em grosso modo}, $\lambda_v(P-Q) \geq (1-\epsilon) \min(\lambda_v(P),\lambda_v(Q))$, ce
qui conduit à $\hat{h}(P-Q) \geq (1-\epsilon) \min(\hat{h}(P),\hat{h}(Q))$ pour $P$, $Q$
entiers) \cite[Lem.\ 3.1--3.3]{MR2220098}. L'avantage de cette approche est
une explicitation de la dépendance en $C$ des bornes; 
nous avons essayé de minimiser cette
dépendance avec soin.

Puis nous avons remarqué qu'on peut forcer $\hat{h}(P-Q)$ (et donc la répulsion) à
augmenter en forçant $P$, $Q$ à être proches l'un de l'autre dans
 $\mathbb{Q}_v^2$ pour
une ou plusieurs places $v$ \cite[Lem.\ 3.1--3.3]{MR2220098}.
On obtient donc une borne (via \cite{MR0514023}) qui dépend d'un paramètre
$t$ qu'on a le droit de fixer -- ce qui correspond à la taille de $p$, si on travaille avec un $v=p$ fini. Dans notre énoncé
(\cite[Thm.\ 3.8]{MR2220098}), $t=0$ correspond à une borne ``pure'' (sans l'utilisation d'un $v$
ou $p$) et $t=1$ correspond à l'utilisation d'un $p$ exactement de la
taille
nécessaire pour donner $90^\circ$ de répulsion.
Si $t=1$, l'exposant de la borne obtenue est celui de \textsc{Bombieri}-\textsc{Pila}, 
à savoir, $1/d$. Nous montrons après que, si on diminue $p$ 
légèrement ($t=1-\delta$), la borne est améliorée. Donc, on arrive toujours à obtenir un exposant meilleur que 
celui dans \cite{MR1016893}.\footnote{Néanmoins, les deux bornes ne sont pas directement 
comparables: \textsc{Venkatesh} et moi travaillons avec la hauteur canonique $\hat{h}$,
\textsc{Bombieri}-\textsc{Pila} avec la hauteur naïve $h_x$. Eliminer cette différence est très loin d'être trivial. Très récemment, mon étudiant en thèse (\textsc{D. Mendes da Costa}) a réussi à l'éliminer, et a donné en conséquence des bornes  supérieures a celles en
\cite{MR1016893} dans le sens stricte (pour les courbes elliptiques).}

\subsection{Applications}

La borne avec $t=0$ peut elle aussi être très utile, surtout pour borner le
nombre total de points dans une courbe elliptique (sans restrictions sur la
hauteur). La dépendance des bornes 
de la hauteur est si légère pour $t=0$ que, en combinant nos méthodes avec des bornes
exponentielles sur la hauteur (\textsc{Hajdu} et \textsc{Herendi} \cite{MR1615334}, dans la tradition
initiée par \textsc{Baker}), nous obtenons le résultat suivant.
\begin{prop}$($\cite[Cor.\ 3.12]{MR2220098}$)$ Soit $E$ une courbe elliptique 
définie sur $\mathbb{Q}$ par une équation de Weierstrass avec coefficients entiers. Alors, le
nombre de points entiers sur $E(\mathbb{Q})$ est
\[O\left(|\Delta|^{0.20070\dotsc}\right),\] où $\Delta$ est le discriminant de l'equation de Weierstrass.
\end{prop}

Il est possible de prouver une borne similaire sur un corps de nombres $K$ arbitraire,
ou bien aussi pour le nombre de points $S$-entiers. (Un point $S$-entier est un point 
$(x,y)\in K^2$ tel que $x$, $y$ peuvent être écrits comme des fractions dont les dénominateurs
sont divisibles seulement par des premiers dans l'ensemble $S$.) Nous en dérivons des bornes sur le nombre de courbes elliptiques de
conducteur donné; la famille de ces courbes est paramétrée par des points $S$-entiers sur une courbe
elliptique.

\begin{prop}\label{prop:loro}$($\cite[Thm. 4.5]{MR2220098}$)$ Le nombre de courbes elliptiques sur $\mathbb{Q}$
de conducteur $N$ est
\begin{equation}\label{eq:crost}
O\left(N^{0.22377\dotsc}\right)\end{equation}
\end{prop}
 
 Nous obtenons aussi (cette fois avec $t\ne 0$ soigneusement optimisé) des
 bornes sur la $3$-torsion $H_3(\mathbb{Q}(\sqrt{D}))$ (le sous-groupe qui
 consiste en les éléments $g$ avec $g^3=1$) dans le groupe de classe de toute extension quadratique $\mathbb{Q}(\sqrt{D})$ de $\mathbb{Q}$. 
  Il y a dans cette situation un cercle vertueux: les bornes sur la 
 $3$-torsion conduisent à des meilleures bornes sur le nombre de points sur une courbe 
 $y^2 = x^3+d$, ce qui conduit à des meilleures bornes sur la $3$-torsion, etc.
 
 \begin{prop}$($\cite[Thm. 4.5]{MR2220098}$)$ Pour tout entier $D$ non nul,
 \begin{equation}\label{eq:procul}
 |H_3(\mathbb{Q}(\sqrt{D}))| \ll |D|^{0.44178\dotsc}\end{equation}
 \end{prop}
 
\textsc{Ellenberg} et \textsc{Venkatesh} \cite{MR2331900} 
ont déjà amélioré (\ref{eq:procul}), ce qui améliore
immédiatement (\ref{eq:crost}), car la preuve de Prop.\ \ref{prop:loro} utilise (\ref{eq:procul}).
 
 \begin{center}
 * * *
 \end{center}

Nos méthodes pour des courbes elliptiques sont aussi 
valides (avec des angles plus grands, donc meilleures
bornes) pour des courbes de genre $\geq 2$. Ce fait est utilisé en
\cite{MR2394896}, qui se sert aussi d'arguments probabilistes pour minimiser
le coût de forcer $P$ et $Q$ à être proches $p$-adiquement (dans le
contexte de courbes superelliptiques $t y^{k-1} = f(x)$, $\deg(f)=k$,
provenant du problème de valeurs polynomiales sans facteurs qui soient
 puissances $(k-1)$-èmes).

\section{La mauvaise distribution force-t-elle l'algébricité?}

Soit $S$ un sous-ensemble de $\mathbb{Z}^n \cap \lbrack
0,N\rbrack^n$. Examinons ce qui se
passe si $S$ est mal distribué modulo $p$ pour beaucoup de nombres premiers $p$.
(Par exemple, soit $S_0\subset \mathbb{Z}$, et supposons que, pour chaque $p$ plus grand
qu'une constante, il y a au moins $0.1 p$ classes de congruence dans lesquelles il n'y a aucun
élément de $S_0$.) Intuitivement, il est clair que les éléments de $S$ sont atypiques, et donc
$S$ doit être petit. Peut-on quantifier et exploiter cette intuition?

Le grand crible et le crible majeur (c'est-à-dire le {\em larger sieve} de \textsc{Gallagher} \cite[Thm. \ 1]{MR0291120})
donnent des bornes pour le nombre d'éléments $|S|$ de l'ensemble $S$ dans ce type de
situation. Par exemple, dans le cas de l'ensemble $S_0$ défini ci-dessus, le crible majeur
nous donne que
\begin{equation}\label{eq:symph8}
|S_0|\ll N^{0.9},\end{equation}
où la constante implicite est absolue. Le grand crible va plus loin dans ce cas: 
\cite[Thm.\ 6]{MR0371840} nous donne
\begin{equation}\label{eq:kindtot}
|S_0|\ll N^{1/2} (\log N)^c, 
\end{equation}
où $c$ et la constante implicite sont absolues.

Ces bornes sont-elles optimales? Il y a des ensembles pour lesquels le
grand crible et le crible majeur donnent tous les deux des bornes tout à fait optimales. Par exemple, soit $S_0$
l'ensemble de tous les carrés jusqu'à $N$. Alors, pour chaque $p$, il y a $(p-1)/2$ 
(donc $>0.1p$) classes
de congruence modulo $p$ dans lesquelles il n'y a aucun élément de $S_0$. Bien entendu,
$|S_0| = \lfloor N^{1/2}\rfloor$. Donc, dans ce cas, (\ref{eq:kindtot}) est à peu près optimal.

En même temps, cet exemple pourrait être une exception: les carrés sont loin d'être un ensemble
arbitraire, puisque, évidemment, ils possèdent une structure algébrique très forte. 
Il est plausible qu'il existe une dichotomie claire \cite[\S 4.2, "Guess"]{MR2508647}:
peut-être,
si un ensemble évite beaucoup de classes de congruence modulo $p$ pour beaucoup de
premiers $p$, soit l'ensemble est très petit, soit il a une forte
corrélation avec une structure algébrique.
\textsc{Venkatesh} et moi avons prouvé un tel résultat en dimension $n=2$:

\begin{theo}
Soit $S\subset \mathbb{Z}^2 \cap \lbrack 0,N\rbrack^2$, $N\geq 1$. Supposons que le nombre
de résidus $\{(x,y) \mo\; $p$ : (x,y)\in S\}$ est $\leq \alpha p$ pour $\alpha>0$ donné et tout $p$.

Alors, pour tout $\epsilon>0$, soit
\begin{itemize}
\item $|S| \ll_{\alpha,\epsilon} N^\epsilon$, soit
\item il y a une courbe algébrique $C$ dans le plan avec degré $O_{\alpha,\epsilon}(1)$
telle que au moins $(1-\epsilon) |S|$ points de $S$ sont sur $C$.
\end{itemize}
\end{theo}
La preuve repose en partie sur de l'idée du crible majeur et en partie sur des méthodes reliées
à celles de \textsc{Bombieri}-\textsc{Pila} \cite{MR1016893}.

\textsc{M. Walsh} a réussi à généraliser le théorème a toutes les dimensions $n>2$ \cite{Walsh}. Le
cas $n=1$ reste ouvert.

\chapter{Théorie analytique et probabiliste des nombres}

Un entier $n$ est dit {\em sans facteurs carrés} si $d^2\nmid n$ pour tout entier $d>1$. Soit
$f\in \mathbb{Z}\lbrack x\rbrack$ un polynôme. Dans la théorie des nombres, il est assez 
souvent convenable ou nécessaire de travailler avec des entiers sans
facteurs
 carrés.
 Y a-t-il une infinité d'entiers $n$ tels que $f(n)$ soit sans facteurs carrés?

Il est évident qu'il est nécessaire de demander comme condition que 
$f(n)\not\equiv 0 \mod p^2$ aie des solutions dans
$\mathbb{Z}/p^2\mathbb{Z}$ pour tout $p$ (ce qui est certainement vrai pour $p$ plus
grand qu'une constante). Il est aussi évidemment nécessaire que $f(n)$ n'aie pas de racines
multiples (donc facteurs répétés, comme $(x+5)^2$). Sous ces conditions, 
il est conjecturé qu'il y a une infinité
d'entiers $n$ tels que $f(n)$ soit sans facteurs carrés. Ceci est trivial pour $\deg(f)=1$, classique pour
$\deg(f)=2$ \cite{MR1512732}, connu pour $\deg(f)=3$ 
(\textsc{Erd\H{o}s}, \cite{MR0056635}) et ouvert
pour $\deg(f)>3$ (mais prouvé dans \cite{MR1654759} 
de façon conditionnelle, sous la conjecture ABC).

La preuve d'\textsc{Erd\H{o}s} pour $\deg(f)=3$ évite astucieusement
plusieurs questions
difficiles en géometrie diophantienne. Peut-être en consequence, \textsc{Erd\H{o}s} a demandé si
on pouvait prouver que $f(p)$ est sans facteurs carrés pour une infinité de premiers $p$: ceci
rend ses propres astuces inutiles, et au moins quelques questions diophantiennes doivent être alors adressées.

\begin{theo}$($\cite{MR2437966}$)$\label{theo:hongo}
Soit $f\in \mathbb{Z}\lbrack x\rbrack$ un polynôme cubique sans racines multiples. Supposons que,
pour tout nombre premier $q$, il y une solution $x\in (\mathbb{Z}/p^2\mathbb{Z})^*$ à
 $f(x)\not\equiv 0 \mo p^2$.
 Alors il y a une infinité de nombres premiers tels que $f(p)$ est sans facteurs carrés.
\end{theo}
En plus, \cite{MR2437966} calcule la proportion (positive, sous les
conditions) de nombres premiers $p$ tels que $f(p)$ soit sans facteurs de carrés.

La partie difficile du problème est de prouver que $f(p)$ peut avoir un grand diviseur
carré du type $q^2$, $q$ premier, seulement pour une proportion $o(1)$ des nombres premiers
$p$. L'expression
\begin{equation}\label{eq:milhaud}
d q^2 = f(p)
\end{equation}
rend le lien avec les courbes elliptiques $d y^2 = f(x)$ évident.
(Dans cette perspective, il est possible de réinterpreter le travail de
\textsc{Hooley} \cite[Ch. IV]{MR0404173} comme une approche exclusivement
locale
à un problème de comptage, i.e., une approche qui ne prend pas avantage
de la non-validité du principle de \textsc{Hasse}.)

Par des méthodes basées sur la répulsion (developpées dans \cite{MR2508647}; voir \S \ref{sec:courbes} de ce mémoire), le nombre de points entiers sur $ d y^2 = f(x)$ est petit
($\ll (\log N)^{O(\epsilon)}$) pour tout $d$ ayant (a) $O(\log \log N)$ facteurs premiers et
(b) un grand diviseur $d_0 | d$ ($d_0\gg N^{1-\epsilon}$) avec $O(\epsilon \log \log N)$
facteurs premiers. Il n'est pas difficile à prouver que ceci est le cas pour une proportion
$1-o(1)$ de tous les $d$, et aussi pour tout $d$ venant de $d = f(p)/q^2$ pour une proportion $1-o(1)$ de toutes les valeurs $p$. On peut donc travailler seulement avec des valeurs de $d$
pour lesquelles $d y^2 = f(x)$ a $(\log N)^{O(\epsilon)}$ solutions. Ceci n'est pas assez pour
prouver le Théo. \ref{theo:hongo}: on doit considérer $\sim N (\log N)^\epsilon$
valeurs de $d$, et on veut une borne totale de grandeur $o(N/\log N)$ -- donc même une borne de $O(1)$
points par courbe ne serait pas assez.

Il est alors nécessaire d'éliminer la majorité des valeurs $d$. Pour faire ceci, on montre que, pour 
$p$ typique, $f(p)$ est dans la queue de la distribution du nombre de diviseurs premiers d'un 
entier, classifiés par leur type de décomposition en $\mathbb{Q}(\alpha)/\mathbb{Q}$,
où $\alpha$ est une racine de $f(\alpha)=0$. C'est-à-dire: $f(p)$ n'est divisible par aucun 
nombre premier inerte dans $\mathbb{Q}(\alpha)/\mathbb{Q}$, mais, par contre, 
en moyenne, il est divisible par
$3$ fois plus de nombres premiers totalement décomposés qu'un entier
moyen. Comme $f(p)$
a ces propriétés atypiques, $d = f(p)/q^2$ ($q$ premier) les a aussi, donc $d$ est atypique. Ceci est équivalent
à dire que $d$ est dans un ensemble peu nombreux; cette équivalence se quantifie de manière précise avec des bornes de {\em grandes déviations}.

Si le groupe de Galois de $\mathbb{Q}(\alpha)/\mathbb{Q}$ est 
le groupe alterné $\Alt(3)$, ces bornes sont
suffisamment fortes pour finir la preuve du Théo. \ref{theo:hongo}.
Ceci est le résultat principal de \cite{MR2394896}, qui résout aussi le
problème analogue pour les valeurs de $f(p)$ ($\deg(f)=k$) sans facteurs
qui
soient
puissances
$(k-1)$-èmes, sous l'hypothèse que le groupe de Galois de 
$\mathbb{Q}(\alpha)/\mathbb{Q}$ soit de {\em haute entropie}
(\cite[Thm. 1.1]{MR2394896}).

Si, par contre, le groupe de Galois de $\mathbb{Q}(\alpha)/\mathbb{Q}$ est 
le groupe symétrique $\Sym(3)$, 
il faut encore montrer
que $d$ est atypique d'une autre façon 
supplémentaire -- la suivante. Le nombre de points entiers sur
$y^2 = f(x) \mod p$ est $p + O(\sqrt{p})$ (\textsc{Hasse}(-\textsc{Weil})), pas exactement $p$; 
la modularité des courbes elliptiques implique que, en moyenne, le terme d'erreur 
$O(\sqrt{p})$ est bien présent dans la réalité. Alors, dans la moyenne pour $p$,  la valeur
$f(n)$ ($n$ entier ou $n$ premier) a un léger biais en faveur ou à
l'encontre d'être un carré modulo $p$.  Ces biais s'accumulent, et
montrent, à nouveau, que, pour $n$ un premier typique, $f(n)$ est
atypique.  Ceci est assez pour réduire notre attention à assez peu de valeurs de $d$, et le 
Théorème \ref{theo:hongo}  est prouvé.

\chapter{Combinatoire additive classique}

\section{Progressions arithmétiques dans les nombres premiers}
Tout ensemble de densité positive dans les entiers contient un nombre infini de progressions
arithmétiques de longueur 3 (théorème de \textsc{Roth}
\cite{MR0051853}). Les nombres premiers contiennent aussi un nombre infini
de progressions arithmétiques de longueur 3 (\textsc{van der Corput}, basé
sur \textsc{Vinogradov}). Il s'agit de deux applications de l'analyse de
Fourier assez différentes. \textsc{B. Green} \cite{MR2180408}
 a prouvé un résultat hybride: tout ensemble de
densité positive dans les nombres premiers contient un nombre infini de progressions
arithmétiques de longueur 3.  (Ceci précède son travail avec Tao, où
 le même résultat est prouvé pour des progressions de longueur arbitraire \cite[Thm. 1.2]{MR2415379}.)

Affaiblir l'hypothèse de densité positive est une tâche difficile pour laquelle le
progrès est graduel: voir \cite{MR889362}, \cite{MR1100788}, \cite{MR1726234}, \cite{MR2403433}, \cite{Sandlog} (vers la conjecture de \textsc{Erd\H{o}s}-\textsc{Tur\'an}) dans le cas de sous-ensembles des entiers.
Pour le cas de sous-ensembles $A$ de l'ensemble de nombres premiers $P$, \textsc{Green}
a prouvé son résultat sous l'hypothèse $\delta_P(N)\gg 
\sqrt{(\log_5 n)/(\log_4 n)}$, où $\delta_P$ est la {\em densité relative}
$\delta_P(N) = |A\cap \lbrack 1,N\rbrack|/|P\cap \lbrack 1, N\rbrack|$ et
$\log_k$ est le logarithme réitéré $k$ fois. 

\textsc{A. de Roton} et moi avons prouvé le résultat suivant.
\begin{theo}[\cite{MR2773330}]
Soit $A\subset P$ avec densité relative \[\delta_P(N)\geq 
\frac{C \log \log \log N}{(\log \log N)^{1/3}},\]
où $C$ est une constante absolue. Alors $A$ a une infinité de progressions arithmétiques
de longueur $3$.
\end{theo}
Il s'agit donc d'une amélioration du résultat de \textsc{Green} par deux logarithmes
($2$ en place de $4$). 

Un des deux logarithmes gagnés vient de l'utilisation de l'{\em enveloping sieve}, une manière d'utiliser des cribles majorants, ainsi nommée dans \cite{MR0404173}; voir aussi les
travaux de \textsc{Ramaré} \cite{MR1375315}, \cite{MR1879673}. Cette tactique a été utilisée à partir de \cite{MR2415379} par \textsc{Green} et \textsc{Tao}, qui ont contribué à la populariser.
 
 Le deuxième log vient d'un aperçu peut-être plus original. Comme il est
 bien connu (au moins depuis \textsc{Roth}), le nombre
 de progressions arithmétiques dans un ensemble $S$ peut s'exprimer
de manière naturelle en termes de la transformée de Fourier de la fonction
caractéristique $1_S$ de $S$:
\[
|\{a,d \in \mathbb{Z}: a, a+d, a+2d\in S\}| = \int_0^1
\left(\widehat{1_S}(x)\right)^2 \widehat{1_S}(- 2 x) dx .
\]
Cette expression change peu si on replace $1_S$ par sa convolution
 avec une fonction $\sigma$ telle que $\hat{\sigma}(\alpha)
 \sim 1$ quand $\widehat{1_S}(\alpha)$ est grand (car alors
$\widehat{1_S\ast\sigma}(\alpha) = \widehat{1_S}(\alpha) \sigma(\alpha)
\sim
 \widehat{1_S}(\alpha)$ quand $\widehat{1_S}(\alpha)$ est grand). Jusqu'à
ici, tout est en essence comme dans \cite{MR2180408}.  
 De là, par contre, il s'avère plus avantageux
 de travailler dans l'espace physique (plutôt que dans l'espace de Fourier, comme dans \cite{MR2180408}). Une convolution par $\sigma$ (ou par toute autre fonction pas trop grande
 en norme $\ell_\infty$) ramollit dans le sens $\ell_2$ l'ensemble $P$ des premiers, grâce à
 des bornes supérieures (fournies par des cribles) sur les auto-corrélations de la fonction 
 caractéristique de $P$. Ceci permet d'extraire de $1_S\ast \sigma$ un ensemble de
 grande densité en $\mathbb{Z}$, auquel il est possible d'appliquer le théorème de 
 \textsc{Roth}, ce qui complète la preuve.
 
\section{Sommes, produits, incidence}\label{subs:soprod}
Le théorème de \textsc{Szemerédi}-\textsc{Trotter} \cite{MR729791}
 borne le nombre d'éléments de l'ensemble $I$
d'incidences entre un ensemble $P$ de points et un ensemble $L$ de droites dans le
plan réel $\mathbb{R}^2$. (Une
{\em incidence} est une paire $(p,\ell)$, où $p\in P$, $\ell \in L$ et $p$ est sur $\ell$.) Par
exemple, dans le cas particulier $|L|\sim |P| \sim N$, la borne $|I| \leq N^2$ est triviale,
la borne $|I|\ll N^{3/2}$ est relativement facile, et la borne $|I|\ll N^{4/3}$ fournie par
\textsc{Szemerédi}-\textsc{Trotter} est en général optimale.

Il est possible de déduire de  \cite{MR729791} (a) une borne inférieure sur
le nombre de lignes déterminées par un ensemble de $N$ points (théorème de
\textsc{Beck}) et (b) un 
théorème de type \textsc{Erd\H{o}s}-\textsc{Szemerédi}, c'est-à-dire, \[\min(|A\cdot A|,|A+A|)\geq |A|^{1+\delta}\]
pour tout ensemble fini $A\subset \mathbb{Z}$ ($\delta =1/4$, \textsc{Elekes} \cite{MR1472816}).

Dans le corps fini $\mathbb{F}_p$, 
\textsc{Bourgain}, \textsc{Katz} et \textsc{Tao} \cite{MR2053599} ont prouvé
une borne de type \textsc{Erd\H{o}s}-\textsc{Szemerédi}:
\begin{equation}\label{eq:hust}
\min(|A\cdot A|,|A+A|)\geq |A|^{1+\delta}\end{equation}
pour tout ensemble fini $A\subset \mathbb{F}_p$ avec $|A|<p^{1-\epsilon}$
 et une constante $\delta>0$ qui dépend seulement de $\epsilon>0$
("théorème somme-produit"; ceci comprend la contribution de \textsc{Konyagin} \cite{Koarx}). Ils arrivent à en déduire une borne de type
\textsc{Szemerédi}-\textsc{Trotter} ($|I|\leq  N^{3/2-\delta'}$) dans $\mathbb{F}_p^2$; malheureusement, la valeur de $\delta'>0$
qui résulte n'a pas été explicitée, et devrait être minuscule.

En partant d'un résultat de type (\ref{eq:hust}) avec $\delta>0$ explicite
(\textsc{Li} \cite{MR2771659}; voir aussi \cite{MR2461864} et al.), \textsc{Rudnev} et moi
avons prouvé \cite{MR2764161} une version du théorème de Beck dans $\mathbb{F}_p^2$
avec un exposant explicite. La preuve repose sur le fait qu'un nombre trop
élevé de triples de points
colinéaires impliquerait un grand nombre de relations linéaires, menant
après un certain effort à une contradiction avec (\ref{eq:hust}). Nous avons déduit de ceci un résultat de type
\textsc{Szemerédi}-\textsc{Trotter} avec $\delta'>0$ explicite 
($\delta' = 1/9278$).

\chapter{Théorie algorithmique des nombres }

Il est possible de vérifier en temps polynomial en $\log n$ si un nombre $n$ est premier.
Ceci est classique (\textsc{Solovay}-\textsc{Strassen}, \textsc{Miller}-\textsc{Rabin}, etc.) si
on permet des algorithmes probabilistes, et connu depuis \cite{MR2123939}
de manière déterministe. Mais comment {\em trouver} un nombre premier --
entre $N$ et $2 N$, par exemple?

Si on permet un algorithme probabiliste, la réponse est évidente: pour trouver un nombre premier
entre $N$ et $2 N$, on n'a qu'à prendre un entier entre $N$ et $2 N$ de manière aléatoire, et
vérifier sa primalité; c' nombre premier de grandeur $\sim N$ est premier avec 
probabilité $\sim 1/\log N$, l'algorithme réussira avec probabilité au moins $1 - \epsilon$ après
$O_\epsilon(\log n)$ itérations. Par contre, si seulement des algorithmes déterministes sont
permis, la question semble devenir très difficile.

Le travail qui suit a été fait dans le contexte d'un projet {\em Polymath}
de collaboration massive. Par une décision du journal, l'article \cite{TCH}
sera publié sous le nom des collaborateurs principaux, à savoir,
\textsc{Tao} (organisateur), \textsc{Croot}, et moi-même.

On pourrait vérifier les entiers $N, N+1, N+2,\dotsc$ l'un après l'autre, mais - sauf si on suppose
la conjecture de \textsc{Cramér} - rien ne garantit l'existence d'un nombre premier dans
$\lbrack N, N+M\rbrack$ pour $M$ petit; même l'hypothèse de Riemann nous donnerait ceci
seulement pour $M>N^{1/2+\epsilon}$.

Si on pouvait déterminer rapidement si un intervalle $\lbrack N,N+M\rbrack$ contient des nombres
premiers, on pourrait trouver facilement un nombre premier dans $\lbrack N, 2N\rbrack$. (``J'ai
un nombre entre $1$ et $100$.'' ``Est-ce qu'il est plus grand que $50$?'', etc.) Si le 
nombre de nombres premiers entre $N$ et $2N$ est impair, il suffit d'avoir un algorithme 
qui détermine si l'intervalle $\lbrack N, N+M\rbrack$ contient un nombre pair ou impair de nombres
premiers. Or
\[2\cdot |\{N<p\leq N+M: \text{$p$ premier}\}| = \mathop{\sum_{N< n \leq N+M}}_{\text{$n$ sans diviseurs carrés}}
 \tau(n) \mod 4,\]
où $\tau(n)$ est le nombre de diviseurs de $n$.
Le problème se réduit donc
à calculer $\sum \tau(n) \mod 4$, ce qu'on peut faire en calculant $\sum \tau(n)$ de manière
complètement précise.

Le {\em problème de Dirichlet} consiste à donner une expression algébrique de
\[\sum_{n\leq N} \tau(n)\] avec un terme d'erreur (inévitablement plutôt grand) à minimiser.
\textsc{Gauss} a prouvé $O(N^{1/2})$, \textsc{Voronoi} (1903) a prouvé
$O(N^{1/3})$, \dots .
Notre problème consiste à donner un algorithme qui calcule $\sum_{n\leq N} \tau(n)$ sans
terme d'erreur, et à minimiser le temps que l'algorithme (malheureusement plutôt lent!) prend.
Nous verrons maintenant que non seulement les idées de \textsc{Gauss} mais aussi celles 
de \textsc{Voronoi} et ses contemporains (notamment \textsc{Vinogradov}) sont utiles pour
cette tâche. 

Comme dans le problème de Dirichlet traditionnel, nous commençons par
\[\sum_{n\leq N} \tau(n) = \sum_{n\leq N} \sum_{d|n} 1 = \sum_{d\leq N} \left\lfloor \frac{N}{d}
\right\rfloor = \text{aire sous hyperbole $y =N/x$}.
\]
Par symétrie (\textsc{Gauss}),  cette aire est deux fois l'aire sous l'hyperbole $y=N/x$ et à gauche
de la ligne $x=\sqrt{N}$, moins $(\sqrt{N})^2 = N$. Donc il faut calculer 
$\sum_{d\leq \sqrt{N}} \lfloor N/d\rfloor .$ Ceci se fait en temps $O(\sqrt{N})$ de manière triviale.
Comment faire mieux?

Nous coupons $\lbrack 0,\sqrt{N}\rbrack$ en sous-intervalles $I$ de longueur $\Delta$ avec
le coût d'un facteur $\sqrt{N}/\Delta$ en temps. Evidemment, $\sum_{d\in I} \lfloor N/d
\rfloor = \sum_{d\in I} N/d - \sum_{d\in I} \{ N/d\}$ et $\sum_{d\in I} N/d$ peut être calculé
rapidement à la précision nécessaire grâce à \textsc{Euler}-\textsc{Maclaurin}. (On peut se 
permettre une erreur totale jusqu'à $1/2$ dans chaque intervalle $I$: l'erreur disparait quand 
l'approximation de $\sum_{d\in I} \lfloor N/d \rfloor \in \mathbb{Z}$ est arrondie à l'entier le
plus proche.) Il reste à calculer $\sum_{d\in I} \{ N/d\}$.

Pour $\Delta$ assez petit, $d\mapsto N/d$ est à peu près linéaire en $I$. On remplace
$d\mapsto N/d$ par une approximation linéaire, avec une approximation diophantienne à
$(N/d_0)' = -N/d_0^2$ ($d_0\in I$) comme pente. On obtient une simple somme linéaire mod
$q$, sauf pour quelques cas spéciaux ($N/d$ proche de $\mathbb{Z}$) à identifier et calculer
de manière individuelle. L'algorithme est rapide et valide pour $\Delta\ll N^{1/6}$,
donc le temps total est $(N^{1/2}/\Delta) \cdot N^\epsilon = N^{1/3+\epsilon}$.

\begin{theo}$($\cite[Thm. 1.2]{TCH} et \cite[\S 2.1]{TCH}$)$
La parité du nombre de nombres premiers entre $N$ et $2N$ peut être
déterminée en temps $O_\epsilon\left(N^{1/3 +\epsilon}\right)$ pour tout 
$\epsilon>0$. Si le nombre de nombres premiers entre $N$ et $2N$ est
impair, un nombre premier entre $N$ et $2 N$ peut être trouvé de manière
déterministe en temps $O_\epsilon\left(N^{1/3 +\epsilon}\right)$.
\end{theo}

La condition que $|\{N<p\leq 2N\}|$ soit impair
 peut être remplacée par des autres conditions 
similaires -- soit en considérant $p$ dans des progressions arithmétiques,
soit
en utilisant la notion de {\em complexité de circuits} \cite[Thm.\ 1.3]{TCH}.

\chapter{Conclusions et perspectives}

{\bf \S 2.} {\em Algorithmes.} Disons que nous avons prouvé que le
diamètre d'un graphe de Cayley $\Gamma(G,A)$ est petit. Il reste à voir si
ce résultat est ou peut devenir constructif. C'est-à-dire: est-il possible
de trouver, pour tout $g\in G$, des éléments $a_1,a_2,\dotsc,a_\ell \in
A\cup A^{-1}$ tels que $a_1\dotsb a_\ell = g$ et $\ell$ est petit. Par
exemple, un algorithme dans un contexte similaire
(approximations pour $A$ fixe) existe pour $SU(2)$
(\textsc{Solovay}-\textsc{Kitaev}). En
$\SL_2(\mathbb{Z}/p\mathbb{Z})$, par contre, si effectivement il y a un algorithme
probabiliste
(\textsc{Larsen} \cite{MR1976231}), il est valide seulement pour un choix très
particulier de $A$. Tant
\cite{MR2415382} que les travaux qui le suivent sont non-constructifs.
Par contre, il n'est encore clair que \cite{HS} ne puisse devenir
constructif -- des parties de sa preuve sont déjà proches à l'être.

Même un cas plus basique que $\SL_2(\mathbb{Z}/p\mathbb{Z})$ reste ouvert.
Soit $G$ un
sous-groupe de Borel de  $\SL_2(\mathbb{Z}/p\mathbb{Z})$ et $A$ un
ensemble de générateurs. Le diamètre de $\Gamma(G,A)$ est $\ll (\log
|G|)^{O(1)}$ grâce à une estimation de \textsc{Konyagin} (\textsc{Bukh} et
\textsc{Helfgott}, non-publié). Cette méthode est malheureusement
non-constructive. Le problème donc demeure; il se réduit à la question
suivante: comment exprimer un
element $x\in \mathbb{Z}/p\mathbb{Z}$ comme une combinaison linéaire de puissances $1$,
$r$, $r^2$,\dots,
$r^\ell$, $r$ donné, $\ell \ll (\log p)^{O(1)}$, avec coefficients
$\ll (\log p)^{O(1)}$? (Ceci est clair pour $r$
fixe et $p\to \infty$; or le problème consiste à en donner une solution
valide pour $r$ arbitraire.)

{\em Diamètre polynomial.} Si l'analogue de la proposition
\ref{prop:key} est faux pour
$G = \Sym(n)$ et $G=\Alt(n)$, la possibilité qu'il existe une borne
polynomiale
$\diam(G) \leq n^{O(1)} = (\log |G|)^{O(1)}$ reste bien réelle. Peuvent
les idées en \cite{HS}
nous mener à ce but? Une fois cette tâche accomplie, il semble plausible
que la conjecture
de \textsc{Babai} pousse être finalement prouvée en général -- en
particulier, pour les groupes
linéaires de rang non-borné.

{\bf \S 3.1} Mon étudiant en thèse (\textsc{D. Mendes da Costa}) est
arrivé à obtenir des bornes
strictement meilleures que \textsc{Bombieri-Pila} pour les courbes
elliptiques sur $\mathbb{Q}$
définies par une forme courte de Weierstrass. Il est encore à voir si une
méthode similaire pourrait
faire le même pour les courbes de genre $>1$.

{\bf \S 3.2} Le résultat principal de \cite{MR2508647} a été généralisé
aux dimensions $n>2$
par \textsc{Walsh} \cite{Walsh}. Le cas $n=1$ reste ouvert. Il est sans
doute très difficile, puisque
une solution fournirait des résultats de loin plus forts que le grand
crible dans la majorité de situations. Du travail exploratoire par
l'auteur et \textsc{J. Wolf} suggère que les (semi-)normes
de \textsc{Gowers} ne sont pas nécessairement le meilleur moyen pour la tâche,
car elles détectent des (graphes de) {\em fonctions} polynomiales plutôt
que des {\em valeurs}
polynomiales.

{\bf \S 4.} \textsc{Heath-Brown} \cite{MR2521494} et \textsc{Browning}
\cite{MR2773215} ont prouvé
l'infinité de valeurs de $f(p)$
($\deg(f)=k$) sans facteurs qui soient puissances $(k-1)$-èmes pour $k\geq
5$. En conséquence,
il reste à prouver seulement le cas $\deg(f)=k=4$, $\Gal_f$ égal à
$\Alt(4)$ ou $\Sym(4)$.

{\bf \S 5.1.}  Des travaux récents par \textsc{Bateman}-\textsc{Katz}
\cite{BatKatz} et \textsc{Sanders} \cite{Sandlog} suggèrent que la conjecture de
\textsc{Erd\H{o}s}-Tur\'an pourrait être prouvée dans le
futur prévisible. Celle-ci impliquerait immédiatement le fait déjà connu
de l'existence d'un nombre infini de progressions arithmétiques de
longueur $3$ dans tout sous-ensemble des nombres premiers de densité
relative positive -- et pourrait donner des bornes meilleures que celles
développées en \cite{MR2180408} et \cite{MR2773330}.

{\bf \S 5.2.} Il y a des théorèmes d'incidence qui peuvent être prouvés
directement dans un plan projectif abstrait, sans référence à un corps de
base. Ceci est exploré en \textsc{Gill}-\textsc{Helfgott} (non-publié). En
particulier, il est possible de prouver une borne d'incidence (de type
assez différent de \textsc{Bateman}-\textsc{Katz}) dans un plan projectif
sans l'axiome de \textsc{Desargues} ni
celui de \textsc{Pappus}. Or dans ce cas il n'y a une structure algébrique
(corps, anneau) qui
corresponde au plan de manière simple et directe. Donc il s'agit d'un
théorème d'incidence qui n'est  pas une paraphrase d'un résultat
algébrique.

{\bf \S 6.} L'article \cite{TCH} a montré que des méthodes combinatoires
développées pour le problème de \textsc{Dirichlet} peuvent être adaptées à
un régime algorithmique. Il reste à voir si des méthodes basées sur les
sommes exponentielles peuvent être adaptées au même régime. (De telles
méthodes donnent les meilleures bornes connues pour les problème de
\textsc{Dirichlet}.)

\bibliographystyle{alpha}
\bibliography{hdr}
\end{document}